	\documentstyle{article}
	\newcommand{\al}{\alpha}
	\newcommand{\del}{\delta}
	\newcommand{\eps}{\epsilon}
	\newcommand{\lam}{\lambda}
	
	\newcommand{\sig}{\sigma}
	\newcommand{\th}{\theta}
	\newcommand{\om}{\omega}
	\newcommand{\Del}{{\it \Delta}}
	\newcommand{\Gam}{{\it \Gamma}}
	\newcommand{\Lam}{{\it \Lambda}}
	
	\newcommand{\Om}{{\it \Omega}}

	\font\bbb=msbm7 
	\font\BBB=msbm10 
	\newcommand{\NN}{{\mbox{\BBB{N}}}}
	\newcommand{\ZZ}{{\mbox{\BBB{Z}}}}
	\newcommand{\zz}{{\mbox{\bbb{Z}}}}

	\newcommand{\RE}{{\mbox{\BBB{R}}}}
	\newcommand{\co}{{\mbox{\bbb{C}}}}
	\newcommand{\CO}{{\mbox{\BBB{C}}}}
	\font\frak=eufm10 at 11 pt
	
	\newcommand{\OO}{{\cal O}}

	\newcounter{sect}
	\newcommand{\sect}[1]{\vspace{4ex}\addtocounter{sect}{1}
		\begin{flushleft}
		{{\large\bf \arabic{sect}. {#1}}}
		\end{flushleft}
		\setcounter{thm}{0}
		\setcounter{equation}{0}
		\def\theequation{\arabic{sect}.\arabic{equation}}}
	
	\newtheorem{thm}{Theorem}[sect]
	\newtheorem{prop}[thm]{Proposition}
	\newtheorem{lemma}[thm]{Lemma}
	\newtheorem{cor}[thm]{Corollary}
	\newtheorem{defn}[thm]{Definition}
	\newtheorem{rmk}[thm]{Remark}
	\newcommand{\proof}[1]{\noindent {\em Proof.}$\quad$ {#1} $\hfill\Box$
				\vspace{2ex}}
	\newcommand{\be}{\begin{equation}}
	\newcommand{\ee}{\end{equation}}
	\newcommand{\bea}{\begin{eqnarray}}
	\newcommand{\eea}{\end{eqnarray}}
	\newcommand{\nno}{\nonumber \\}
	\newcommand{\sep}[1]{\!\!\!\! &{#1}& \!\!\!\! }
	\newcommand{\eq}{\sep{=}}
	\newcommand{\vc}{\sep{ }}

	\newcommand{\m}{\backslash}
	\newcommand{\inv}[1]{\frac{1}{#1}}
	\newcommand{\hf}{{\textstyle\inv{2}}}
	\newcommand{\e}[1]{e^{{#1}}}
	\newcommand{\ii}{\sqrt{-1}}
	\newcommand{\dr}{d}
	
	\newcommand{\tr}{\,{\rm tr}}
	\newcommand{\Tr}{\,{\rm Tr}}

	\newcommand{\set}[2]{\{{#1}\,|\,{#2}\}}

	\newcommand{\kb}{{\bar{k}}}
	\newcommand{\lb}{{\bar{l}}}
	\newcommand{\klb}{_{k\lb}}
	\newcommand{\zb}{{\bar{z}}}

	\newcommand{\lpk}{{\lam^p_k}}
	\newcommand{\prodk}{\prod^n_{k=1}}
	\newcommand{\sumk}{\sum^n_{k=0}}
	\newcommand{\sump}{\sum_{p\in F}}
	\newcommand{\Omk}{{\Om^{0,k}}}
	\newcommand{\thT}{\th+\ii T}
	\newcommand{\two}[4]{\left\{    \begin{array}{ll}
					{#1}, & {\mbox{if }} {#2}, \\
					{#3}, & {\mbox{if }} {#4}
					\end{array}     \right.}
	\newcommand{\ka}{K\"ahler }
	\newcommand{\hme}{{H^k(M,\OO(E))}}
	\newcommand{\BTp}{{\cal B}^{2,p}_{T^2}}
	\newcommand{\CTp}{{\cal C}^{2,p}_{T^2}}
	\newcommand{\CT}{{\cal C}^2_{T^2}}
	\newcommand{\HT}{{\cal H}^2_{T^2}}
	\newcommand{\QTp}{Q^p_{T^2}}
	\newcommand{\Rinv}{R^{-1}_p(\th)}
	\newcommand{\pdr}{\partial}
	\newcommand{\hfal}{\frac{\al}{2}}
	\newcommand{\pdrb}{\bar{\pdr}}
	\newcommand{\Boxb}{\stackrel{-\!\!\!\!-}{\Box}}
	\newcommand{\form}{\dr J\dr h}
	\newcommand{\Lhat}{\hat{L}_V}
	\newcommand{\Vhat}{\hat{V}}
	\newcommand{\con}{S_2({\textstyle-\frac{\pi}{2}})}
	\newcommand{\uT}{_{u,T}}
	\newcommand{\uTu}{_{u,T/u}}
	\newcommand{\Thu}{_{Th/u^2}}

\begin{document}
$\!\,{}$

        \vspace{-5ex}

        \begin{flushright}
{\tt dg-ga/9602007} (February, 1996)\\
ICTP preprint IC/96/29
        \end{flushright}

	\begin{center}
	{\LARGE\bf Equivariant Holomorphic Morse Inequalities I:\\
	\vspace{1ex}
		A Heat Kernel Proof}\\

	\vspace{4ex}
	{\large\rm Varghese Mathai}\\
	{\em Department of Pure Mathematics, The University of Adelaide,
	Adelaide, South Australia, Australia}\\

	\vspace{2ex}

	{\large\rm Siye Wu}\footnote{Current address: Mathematical Sciences
	Research Institute, 1000 Centennial Drive, Berkeley, CA 94720, USA} \\
	{\em Mathematics Section, International Centre for Theoretical Physics,
	Trieste I-34100, Italy}
	\end{center}

	\vspace{3ex}

	\begin{quote}
{\small {\bf Abstract.}
Assume that the circle group acts holomorphically on a compact K\"ahler 
manifold with isolated fixed points and that the action can be lifted
holomorphically to a holomorphic vector bundle.
We give a heat kernel proof of the equivariant holomorphic Morse inequalities.
We use some techniques developed by Bismut and Lebeau.
These inequalities, first obtained by Witten using a different argument,
produce bounds on the multiplicities of weights occurring in the twisted
Dolbeault cohomologies in terms of the data of the fixed points.}

	\end{quote}

	\vspace{2ex}

\sect{Introduction}

Morse theory obtains topological information of manifolds from
the critical points of the functions.
Let $h$ be a Morse function on a compact manifold of real dimension $n$
and suppose that $h$ has isolated critical points only.
Let $m_k$ ($0\le k\le n$) be the $k$-th Morse number, the number of
critical points of Morse index $k$.
The Lefschetz fixed-point formula says that the alternating sum of $m_k$
is equal to that of the Betti numbers $b_k$:
	\be
\sumk(-1)^km_k=\sumk(-1)^kb_k.
	\ee
Replacing $(-1)$ by $t$, we get two polynomials (Morse and Poincar\'e
polynomials, respectively) in $t$ that are equal at $t=-1$, i.e.,
	\be\label{strong}
\sumk m_kt^k=\sumk b_kt^k+(1+t)q(t)
	\ee
for some polynomial $q(t)=\sumk q_kt^k$.
The (strong) Morse inequalities assert that $q(t)\ge0$ in the sense that
$q_k\ge0$ for every $0\le k\le n$.

In a celebrated paper [\ref{W0}], Witten showed that the cohomology groups
of the de Rham complex $(\Om^*,\dr)$ can be viewed as the space of ground
states of a supersymmetric quantum system and that the Morse inequalities
can be derived by using a deformation
	\be
\dr_h=\e{-h}\dr\e{h}
	\ee
which preserves the supersymmetry.
The idea was used by Bismut [\ref{Bi}] to give a heat kernel proof of the
Morse inequalities.
Let $\dr^* $ and $\dr^*_h$ be the (formal) adjoints of $\dr$ and $\dr_h$
(with a choice of Riemannian metrics), respectively, and let
	\be
\Del=\{\dr,\dr^*\}\quad\mbox{and}\quad\Del_h=\{\dr_h,\dr^*_h\}
	\ee
be the corresponding Laplacians.
(We adopt the standard notations of operator (anti-)commutators
$\{A,B\}=AB+BA$ and $[A,B]=AB-BA$.)
By Hodge theory,
	\be
\sumk(-1)^k\Tr_{\Om^k}\exp(-u^2\Del)=\sumk(-1)^kb_k
	\ee
for any $u>0$; this is in fact the starting point of the heat kernel proof
of the index theorem.
Similarly, after replacing $(-1)$ by $t$, we obtain
	\be\label{heatker}
\sumk t^k\Tr_{\Om^k}\exp(-u^2\Del)=\sumk b_kt^k+(1+t)q_u(t).
	\ee
It is a straightforward consequence of Hodge theory that the polynomial
$q_u(t)\ge0$.
(See for example [\ref{Bi}, Theorem 1.3]. A slightly different method is
used to show the equivariant version in Lemma \ref{MS} below.)
Since $(\Om^*,\dr_h)$ defines the same cohomology groups as $(\Om^*,\dr)$,
we can replace the heat kernels in (\ref{heatker}) by those associated to
the deformed Laplacian $\Del_h$.
It turns out that
	\be
\lim_{T\to+\infty}\lim_{u\to+0}\Tr_{\Om^k}\exp(-u^2\Del\Thu)=m_k
\quad(0\le k\le n);
	\ee
the (strong) Morse inequalities then follow.
The heart of the proof is that as $u\to0$ the heat kernel is localized near
the critical points of $h$, around which the operator consists of $n$ copies
of (supersymmetric) harmonic oscillators whose heat kernels are given by
Mehler's formula.

Witten [\ref{W}] also introduced a holomorphic analog of [\ref{W0}].
Let $M$ be a compact \ka manifold of complex dimension $n$ and let $E$ be
a holomorphic vector bundle over $M$.
Let $\hme$ be the cohomology groups with coefficients in the sheaf of
holomorphic sections of $E$, calculated from the twisted Dolbeault complex 
$(\Om^{0,*}(M,E),\pdrb_E)$.
Suppose that the circle group $S^1$ acts holomorphically and effectively 
on $M$ preserving the \ka structure and that the action can be lifted
holomorphically to $E$.
Then $\e{\ii\th}$ also acts on the space of sections by sending a section
$s$ to $\e{\ii\th}\circ s\circ\e{-\ii\th}$.
The induced action on $\Om^{0,*}(M,E)$ commutes with the operator $\pdrb_E$.
Thus we obtain representations of $S^1$ on $\hme$;
the multiplicities of weights of $S^1$ in each cohomology group will
be the subject of our investigation.
The $S^1$-action on $(M,\om)$ is clearly symplectic: let $V$ be the vector
field on $M$ that generates the $S^1$-action, then $L_V\om=0$.
If the fixed-point set $F$ of $S^1$ on $M$ is non-empty,
then the $S^1$-action is Hamiltonian [\ref{Fr}], i.e.,
there is a moment map $h\colon M\to\RE$ such that $i_V\om=dh$.
We further assume that $F$ contains isolated points only.
It is well-known that all the Morse indices are even and hence by the
lacunary principle, $h$ is a perfect Morse function: $m_{2k-1}=b_{2k-1}(=0)$,
and $m_{2k}=b_{2k}$ ($0\le k\le n$).
However a refined statement is possible because of the complex structure.
For each $p\in F$, $S^1$ acts on $T_pM$ by the isotropic representation;
let $\lam^p_1,\cdots,\lam^p_n\in\ZZ\m\{0\}$ be the weights.
We define the {\em orientation index} $n^p$ of the fixed point $p\in F$ as 
the number of weights $\lam^p_k<0$; the Morse index of $h$ at $p$ is then
$2(n-n_p)$.
(We need to explain our convention in a simple (but non-compact) 
example $M=\CO$, $\om=\frac{\ii}{2}\dr z\wedge\dr\zb$,
with an $S^1$-action of weight $\lam\in\ZZ\m\{0\}$.
Since $V=\ii\lam(z\frac{\pdr}{\pdr z}-\zb\frac{\pdr}{\pdr\zb})$,
we have $h=-\hf\lam|z|^2$.
Also, the weight of the $S^1$-action on the function $z^k$ 
(a section of the trivial bundle) is $-k\lam$ ($k\in\ZZ$, $k\ge0$);
this leads to a sign convention different from [\ref{W}] in the main result.)
Furthermore, $S^1$ acts on the fiber $E_p$ over $p\in F$.
It is useful to recall a notation in [\ref{W}].
If the group $S^1$ has a representation on a finite dimensional complex
vector space $W$, let $W(\th)$ ($\th\in\RE$) be its character.
For example, we denote $E_p(\th)=\tr_{E_p}\e{\ii\th}$ and
$H^k(\th)=\tr_{\hme}\e{\ii\th}$.
The analog of the Lefschetz formula is the fixed-point formula of Atiyah
and Bott [\ref{AB}], which we write as an equality of alternating sums
[\ref{W0}]:
	\be
\sump(-1)^{n_p}E_p(\th)\prod_{\lpk>0}\inv{1-\e{-\ii\lpk\th}}\prod_{\lpk<0}
\frac{\e{-\ii|\lpk|\th}}{1-\e{-\ii|\lpk|\th}}=\sumk(-1)^kH^k(\th).
	\ee
It turns out that if $(-1)$ is replaced by $t$, the analog of strong Morse
inequalities like (\ref{strong}) holds.
We need the following

\begin{defn}
Let $q(\th)=\sum_{m\in\zz}q_m\e{\ii m\th}\in\RE((\e{\ii\th}))$ be a formal
character of $S^1$, we say $q(\th)\ge 0$ if $q_m\ge 0$ for all $m\in\ZZ$;
let $Q(\th,t)=\sum_{k=0}^nq_k(\th)t^k\in\RE((\e{\ii\th}))[t]$
be a polynomial of degree $n$ with coefficients in $\RE((\e{\ii\th}))$,
we say $Q(\th,t)\ge0$ if $q_k(\th)\ge0$ for all $k$.
For two such polynomials $P(\th,t)$ and $Q(\th,t)$, we say
$P(\th,t)\le Q(\th,t)$ if $Q(\th,t)-P(\th,t)\ge0$.
\end{defn}

Using a holomorphic version of supersymmetric quantum mechanics, Witten 
[\ref{W}] derived the following

\begin{thm}\label{main}
Suppose $M$ is a compact \ka manifold on which $S^1$ acts holomorphically
preserving the \ka form and with non-empty and discrete fixed-point set
and suppose that the $S^1$-action can be lifted holomorphically to a
holomorphic vector bundle $E$ over $M$.
Then under the above assumptions and notations, we have\\
1. Weak equivariant holomorphic Morse inequalities:
	\bea
\vc H^k(\th)\le\sum_{p\in F,n_p=k}E_p(\th)\prod_{\lpk>0}\inv{1-\e{-\ii\lpk\th}}
\prod_{\lpk<0}\frac{\e{-\ii|\lpk|\th}}{1-\e{-\ii|\lpk|\th}},
	\label{weak+}	\\
\vc H^k(\th)\le\sum_{p\in F,n_p=k}E_p(\th)\prod_{\lpk>0}\frac{\e{\ii\lpk\th}}
{1-\e{\ii\lpk\th}}\prod_{\lpk<0}\inv{1-\e{\ii|\lpk|\th}} \quad(0\le k\le n);
	\label{weak-}
	\eea
2. Strong equivariant holomorphic Morse inequalities:
	\bea
\sump t^{n_p}E_p(\th)\prod_{\lpk>0}\inv{1-\e{-\ii\lpk\th}}\prod_{\lpk<0}\frac{
\e{-\ii|\lpk|\th}}{1-\e{-\ii|\lpk|\th}}\eq\sumk t^kH^k(\th)+(1+t)Q^+(\th,t),
	\label{strong+}		\\
\sump t^{n-n_p}E_p(\th)\prod_{\lpk>0}\frac{\e{\ii\lpk\th}}{1-\e{\ii\lpk\th}}
\prod_{\lpk<0}\inv{1-\e{\ii|\lpk|\th}}\eq\sumk t^kH^k(\th)+(1+t)Q^-(\th,t),
	\label{strong-}
	\eea
where $Q^\pm(\th,t)\ge0$;\\
3. Atiyah-Bott fixed-point theorem:
	\be\label{ab}
\sump\frac{E_p(\th)}{\prodk(1-\e{-\ii\lpk\th})}=\sumk(-1)^kH^k(\th).
	\ee
\end{thm}

\proof{Clearly the weak inequalities (\ref{weak+}) and (\ref{weak-}) follow
from the strong ones (\ref{strong+}) and (\ref{strong-}), respectively.
The index formula (\ref{ab}) can be recovered by setting $t=-1$ in either
(\ref{strong+}) or (\ref{strong-}).
Furthermore, we obtain (\ref{strong+}) from (\ref{strong-}) after reversing
the $S^1$-action and replacing $\th$ by $-\th$.
The whole paper is devoted to a heat kernel proof of (\ref{strong-}).}

The cohomology groups $\hme$ as representation spaces of $S^1$
depend only on the ($S^1$-invariant) holomorphic structure on $E$.
We can choose an $S^1$-invariant Hermition form on $E$ and let
$\dr_E=\pdr_E+\pdrb_E$ be the unique compatible holomorphic connection.
To simplify notations, we drop the subscript $E$ but keep in mind that
	\be
\pdr^2=\pdrb^2=0\quad\quad\mbox{and}\quad\quad
\dr^2=\{\pdr,\pdrb\}=\Om\wedge\cdot\,,
	\ee
where the curvature $\Om$ is a $(1,1)$-form on $M$ with values in
${\rm End}(E)$.
Let $\dr^*$, $\pdr^*$, $\pdrb^*$ be the (formal) adjoints of $\dr$, $\pdr$,
$\pdrb$, respectively and let
	\be
\Del=\{\dr,\dr^*\},\quad\Box=\{\pdr,\pdr^*\},\quad\Boxb=\{\pdrb,\pdrb^*\}
	\ee
be the corresponding Laplacians.
Following [\ref{W}], we deform the $\pdrb$ operator and its Laplacian by
	\be
\pdrb_h=\e{-h}\pdrb\e{h},\quad\pdrb^*_h=\e{h}\pdrb^*\e{-h},
\quad\Boxb_h=\{\pdrb_h,\pdrb^*_h\}.
	\ee
The analog of (\ref{heatker}) holds, where $b_k$ should be replaced by
$\dim\hme$ ($0\le k\le n$).
Contrary to the treatment of ordinary Morse theory in [\ref{W0}, \ref{Bi}],
the limit of $\Tr_{\Omk(M,E)}\exp(-u^2\Boxb\Thu)$ as $u\to0$ does not exist.
To see this, we observe that [\ref{W}] (see also formulas (\ref{mag})
and (\ref{w0}) below) up to a (bounded) $0$-th order operator,
$u^2\Boxb\Thu$ is equal to $\hf u^2\Del\Thu+\ii T\Lhat$,
where $\Lhat$ is the infinitesimal action of the circle group $S^1$.
Since $\Lhat$ is an (unbounded) first order differential operator,
the analysis of [\ref{Bi}] that shows localization of heat kernels
does not go through.
{}From the physics point of view, the operator $u^2\Boxb\Thu$ near
a critical point of $h$ is the Hamiltonian operator of a (supersymmetric)
charged particle in a uniform magnetic field.
Therefore the wave functions, and hence the heat kernel, do not localize to
any point no matter how strong the magnetic field is.
However since $\Lhat$ commutes with $u^2\Boxb\Thu$, 
we can restrict the latter to an eigenspace of the former.
The situation changes drastically because $\Lhat$ is then a constant.
Physically, this amounts to fixing the angular momentum with respect to a
given point, which does localize the wave functions to that point in the
strong field limit.
Therefore in this holomorphic setting, we are naturally lead to consider the
equivariant heat kernel and consequently, equivariant Morse-type inequalities.

The inequalities due to Demailly [\ref{D}] have also been referred to
in the literature as holomorphic Morse inequalities.
The important difference with our case is that Demailly's inequalities
do not require a group action and are asymptotic inequalities,
as the tensor power of a holomorphic line bundle gets large,
whereas the inequalities which we consider are for a fixed holomorphic
vector bundle with a holomorphic $S^1$-action, and are not merely asymptotic.

In section 2, we study various deformations of the Laplacians on \ka manifolds.
In particular, the operator $\Boxb_h$ is calculated explicitly.
We also compare two other deformations $\Boxb_v$ and $\Boxb_{\ii v}$,
which are used in studying complex immersions [\ref{BL}] and
holomorphic equivariant cohomology groups [\ref{L}].
Roughly speaking, the operators $\hf\Del_h$, $\Boxb_v$ and $\Boxb_{\ii v}$
form a triplet of a certain $SU(2)$ group.
In section 3, we use the technique of [\ref{BL}] to show that as $u\to0$,
the smooth heat kernel associated to the operator 
$\exp(-u^2\Boxb\Thu+\ii T\Lhat)$ ($u>0$, $T>0$) is localized near the
fixed-point set $F$, and when $F$ is discrete, the equivariant heat kernel
can be approximated by the using the operators with coefficients frozen
at the fixed points.
The result of the previous section is used to relate by a unitary conjugation
the operator $-u^2\Boxb\Thu+\ii T\Lhat$ to $-u^2\Boxb\Thu$ that appears in
[\ref{BL}] (but restricted to a certain subspace) plus a $0$-th order operator
$-\ii Tr_V$ (as $u\to0$) whose action does not depend on the differential
forms.
This has enabled us to follow the analysis of [\ref{BL}] closely, though
a more direct approach without using the conjugation also seems possible.
In section 4, we calculate the equivariant heat kernel of the linearized 
problem using Mehler's formula and then deduce the (strong) equivariant
holomorphic Morse inequalities (\ref{strong-}) by taking the limit
$T\to+\infty$.
Unlike the argument using small eigenvalues [\ref{W}], the $0$-th order
operator $r_V$ plays a crucial role in the heat kernel calculation.

In a separate paper [\ref{Wu}], equivariant holomorphic Morse inequalities
with torus and non-Abelian group actions are established and are applied to
toric and flag manifolds.
The situations with non-isolated fixed points are left for further
investigation.

\sect{Deformed Laplacians on \ka manifolds}

Recall that $E$ is a holomorphic Hermitian vector bundle over a
compact \ka manifold $(M,\om)$.
(The Hermition structure is needed in the proof but not in the statement
of Theorem~\ref{main}.)
Let $\Lam_+=\om\wedge\cdot$ be the exterior  multiplication of $\om$ on
$\Om^{*,*}(M,E)$ and $\Lam_-=\Lam^*_+$, its adjoint.
Then $\Lam_3=\hf[\Lam_+,\Lam_-]$ preserves the bi-grading of $\Om^{*,*}(M,E)$.
In fact, the action of $\Lam_3$ on $\Om^{p,q}(M,E)$ is  $\hf(p+q-n)$,
hence $[\Lam_3,\Lam_\pm]=\pm\Lam_\pm$.
Set $\Lam_1=\hf(\Lam_++\Lam_-)$ and $\Lam_2=-\frac{\ii}{2}(\Lam_+-\Lam_-)$,
then $\Lam_a$ ($a=1,2,3$) satisfy the standard $\mbox{\frak{su}}\,(2)$
commutation relations
	\be
[\Lam_a,\Lam_b]=\ii\eps_{abc}\Lam_c.
	\ee
(See for example [\ref{GH}].)
So there is a unitary representation of $SU(2)$ on $\Om^{*,*}(M,E)$;
let $S_a(\al)=\e{\ii\al\Lam_a}$ be the corresponding group elements.
We now introduce a slightly more generalized setup.

\begin{defn}
Let $\sig\in\Om^{1,1}(M,E)$ be a real-valued $(1,1)$-form.
Set $\Lam_+(\sig)=\sig\wedge\cdot\,$, $\Lam_-(\sig)=\Lam^*_+$,
$\Lam_1(\sig)=\hf(\Lam_+(\sig)+\Lam_-(\sig))$,
$\Lam_2(\sig)=-\frac{\ii}{2}(\Lam_+(\sig)-\Lam_-(\sig))$ and
$\Lam_3(\sig)=\hf[\Lam_+,\Lam_-(\sig)]$($=-\hf[\Lam_-,\Lam_+(\sig)]$).
\end{defn}

\begin{rmk}
{\rm In computations, it is sometimes convenient to introduce local complex
coordinates $\{z^k,k=1,\cdots,n\}$ on $M$.
The \ka form $\om=\om\klb\dr z^k\wedge\dr\zb^\lb$ is related to the metric
$g=g\klb\dr z^k\otimes\dr\zb^\lb$ by $\om\klb=\ii g\klb=-\om_{\lb k}$.
Let $e^k$, $e^\lb$ be the multiplications by $\dr z^k$, $\dr\zb^\lb$,
and $i_k$, $i_\lb$, the contractions by $\frac{\pdr}{\pdr z^k}$,
$\frac{\pdr}{\pdr\zb^\lb}$, respectively.
Clearly they satisfy the following anti-commutation relations:
$\{e^k,i_l\}=\del^k_l$, $\{e^\kb,i_\lb\}=\del^\kb_\lb$ and others $=0$.
If a $(1,1)$-form $\sig=\sig\klb\dr z^k\wedge\dr\zb^\lb$ is real-valued,
then all the $\sig\klb$'s are purely imaginary.
Setting $i^k=g^{k\lb}i_\lb$ and $i^\lb=g^{k\lb}i_k$, we have
$\Lam_+(\sig)=\sig\klb e^ke^\lb$, $\Lam_-(\sig)=-\sig\klb i^ki^\lb$,
$\Lam_1(\sig)=\hf\sig\klb(e^ke^\lb-i^ki^\lb)$,
$\Lam_2(\sig)=-\frac{\ii}{2}\sig\klb(e^ke^\lb+i^ki^\lb)$,
and $\Lam_3(\sig)=-\frac{\ii}{4}\sig\klb([e^k,i^\lb]+[e^\lb,i^k])$.}
\end{rmk}

\begin{lemma}
	\bea
S_1(\al)^{-1}\Lam_3(\sig)S_1(\al)\eq\cos\al\Lam_3(\sig)-\sin\al\Lam_2(\sig),\\
S_2(\al)^{-1}\Lam_3(\sig)S_2(\al)\eq\cos\al\Lam_3(\sig)+\sin\al\Lam_1(\sig),
						\label{conlam}		\\
S_3(\al)^{-1}\Lam_3(\sig)S_3(\al)\eq\Lam_3(\sig).
	\eea
\end{lemma}

\proof{A straight-forward calculation using the above anti-commutation 
relations shows that $[\Lam_a,\Lam_b(\sig)]=\ii\eps_{abc}\Lam_c(\sig)$.
This means that $\{\Lam_a(\sig)\}$ is an $SU(2)$ triplet.
Hence the result.}

It is clear that the Hodge relations (see for example [\ref{GH}])
	\bea
\vc[\Lam_-,\pdr]=\ii\pdrb^*,\quad\quad [\Lam_-,\pdrb]=-\ii\pdr^* \label{com1}\\
\vc[\Lam_+,\pdr^*]=\ii\pdrb,\quad\quad [\Lam_+,\pdrb^*]=-\ii\pdr \label{com2}
	\eea
still hold after coupling to the vector bundle $E$.
Moreover, we have the Bochner-Kodaira-Nakano identities
	\be
\Del=\Box\,+\Boxb,\quad\quad\Boxb-\,\Box=2\Lam_3(\ii\Om),
	\ee
which are consequences of (\ref{com1}), (\ref{com2}) and the graded
Jacobi identities.
(Since $E$ is a holomorphic Hermitian bundle, $\ii\Om$ is a $(1,1)$-form valued
in the subset of ${\rm End}(E)$ which consists of self-adjoint endomorphisms.)
These results have been generalized to non-\ka situations by [\ref{D0}].
When $E$ is a flat bundle, we recover the usual relation $\Box=\Boxb=\hf\Del$.

\begin{lemma}
	\bea
S_1(\al)^{-1}\Boxb S_1(\al)\eq
\Boxb-(1-\cos\al)\Lam_3(\ii\Om)-\sin\al\Lam_2(\ii\Om),\\
S_2(\al)^{-1}\Boxb S_2(\al)\eq
\Boxb-(1-\cos\al)\Lam_3(\ii\Om)+\sin\al\Lam_1(\ii\Om)\label{conbox},\\
S_3(\al)^{-1}\Boxb S_3(\al)\eq\Boxb.
	\eea
\end{lemma}

\proof{From (\ref{com1}) and (\ref{com2}) we deduce that
$S_1(\al)^{-1}\pdrb S_1(\al)=\cos\hfal\pdrb-\sin\hfal\pdr^*$.
Therefore
	\bea
S_1(\al)^{-1}\Boxb S_1(\al)
\eq\{\cos\hfal\pdrb-\sin\hfal\pdr^*,\cos\hfal\pdrb^*-\sin\hfal\pdr\}	\nno
\eq\cos^2\hfal\Boxb+\sin^2\hfal\Box-\cos\hfal\sin\hfal(\{\pdr,\pdrb\}
+\{\pdr^*,\pdrb^*\})							\nno
\eq\Boxb-(1-\cos\al)\Lam_3(\ii\Om)-\sin\al\Lam_2(\ii\Om).
	\eea
The second formula follows in the same fashion from
$S_2(\al)^{-1}\pdrb S_2(\al)=\cos\hfal\pdrb-\ii\sin\hfal\pdr^*$.
The last one is because $\Boxb$ preserves the bi-grading.}


We now equip $M$ with a holomorphic $S^1$-action which preserves the \ka
structure, hence both the complex structure $J$ and the Riemannian metric $g$.
The holomorphic condition $L_VJ=0$ and the Killing equation $L_Vg=0$ read,
in components,
	\be\label{vect}
V_{k;l}=V_{\kb;\lb}\quad\quad\mbox{and}\quad\quad V_{k,\lb}+V_{\lb,k}=0,
	\ee
respectively.
As explained in section 1, we assume that the $S^1$ fixed-point set $F$
is non-empty.
In this case, there is a moment map $h\colon M\to\RE$ satisfying
$i_V\om=\dr h$, or $h_{,k}=-\ii V_k$ and $h_{,\kb}=\ii V_\kb$.
The equations in (\ref{vect}) are equivalent to
	\be\label{fun}
h_{,k;l}=h_{,\kb;\lb}=0\quad\quad\mbox{and}\quad\quad h_{,k,\lb}=h_{,\lb,k}.
	\ee
(The second part is of course the symmetry of the Hessian.)
Also notice the real-valued $(1,1)$-form
	\be
\form=\dr i_Vg=-2\ii h_{,k;\lb}\dr z^k\wedge\dr\zb^\lb.
	\ee

We further assume that the $S^1$-action can be lifted holomorphically
to the bundle $E$.
We can choose an $S^1$-invariant Hermitian form on $E$.
Then the connection $\dr=\dr_E$ is also $S^1$-invariant.
The group element $\e{\ii\th}\in S^1$ acts on a section $s$ by
$s\mapsto\e{\ii\th}\circ s\circ\e{-\ii\th}$.
Let $\Lhat$ be the infinitesimal generator of this $S^1$-action on
$\Om^{*,*}(M,E)$ and let $L_V=\{i_V,\dr\}$ be the standard Lie derivative.
Then the operator
	\be\label{rv}
r_V=\Lhat+L_V
	\ee
is an element of $\Gam(M,{\rm End}(E))$.
Over the fiber of a fixed point $p\in F$, $r_V(p)$ is simply
the representation of ${\rm Lie}(S^1)$ on $E_p$;
this is independent of the choice of the connections on $E$.

\begin{rmk}
{\rm 1. $\Boxb_h$ commutes with the $S^1$-action.
Since the connection, the complex structure, and the moment map $h$
are all $S^1$-invariant, we get $[\Lhat,\dr]=0$, $[\Lhat,\pdrb]=0$ and
$[\Lhat,\pdrb_h]=0$.
Taking the adjoint, we get $[\Lhat,\pdrb^*_h]=0$ and $[\Lhat,\Boxb_h]=0$.\\
2. $\Boxb_h$ also commutes with a $U(1)$ subgroup of $SU(2)$.
Since $\Boxb_h$ preserves the bi-grading, $[\Lam_3,\Boxb_h]=0$,
hence $S_3(\al)^{-1}\Boxb_h S_3(\al)=\Boxb_h$.}
\end{rmk}

\begin{prop}
	\be\label{mag}
\Boxb_h=\Boxb+\hf|\dr h|^2-\Lam_3(\form)-\ii r_V+\ii\Lhat.
	\ee
\end{prop}

\proof{Let $D_k$, $D_\lb$ be the covariant derivative along 
$\frac{\pdr}{\pdr z^k}$, $\frac{\pdr}{\pdr\zb^\lb}$, respectively.
Then $\pdrb=e^\lb D_\lb$, $\pdrb^*=-i^kD_k$ and 
$\pdrb_h=\pdrb+e^\lb h_{,\lb}$, $\pdrb^*_h=\pdrb^*+i^kh_{,k}$.
So
	\bea\label{cal1}
 \Boxb_h\eq\{\pdrb,\pdrb^*\}+\{e^\lb,i^k\}h_{,\lb}h_{,k}
		+\{\pdrb,i^kh_{,k}\}+\{\pdrb^*,e^\lb h_{,\lb}\}		\nno
	\eq\Boxb+g^{k\lb}h_{,k}h_{,\lb}+(h_{,k,\lb}\,e^\lb i^k+h_{,k}D^k)
		-(h_{,\lb,k}i^ke^\lb+h_{,\lb}D^\lb)			\nno
	\eq\Boxb+\hf|\dr h|^2-\Lam_3(\dr J\dr h)+h_{,k,\lb}(e^\lb i^k-e^ki^\lb)
		-\ii(V_kD^k+V_\lb D^\lb).
	\eea
On the other hand,
	\bea\label{cal2}
L_V\eq\{\pdr+\pdrb,V_ki^k+V_\lb i^\lb\}					\nno
   \eq V_kD^k+V_\lb D^\lb+V_{k,\lb}\,e^\lb i^k+V_{\lb,k}e^ki^\lb	\nno
   \eq V_kD^k+V_\lb D^\lb+\ii h_{,k,\lb}(e^\lb i^k-e^ki^\lb).
	\eea
(\ref{mag}) follows from (\ref{cal1}), (\ref{cal2}) and (\ref{rv}).}

We also define two different deformations.
Let $v=V^{1,0}$ be the holomorphic component of $V$. Set
	\be
\pdrb_v=\pdrb+i_v,\quad\quad\Boxb_v=\{\pdrb_v,\pdrb^*_v\}
	\ee
and
	\be
\pdrb_{\ii v}=\pdr+\ii i_v,\quad\quad
\Boxb_{\ii v}=\{\pdrb_{\ii v},\pdrb^*_{\ii v}\}.
	\ee
Then straightforward calculations similar to what leads to (\ref{mag}) yield
	\be
\Boxb_v=\Boxb+\hf|\dr h|^2+\Lam_1(\form)
	\ee
and
	\be
\Boxb_{\ii v}=\Boxb+\hf|\dr h|^2+\Lam_2(\form).
	\ee
It is also interesting to compare the deformation $\Del_h$ in [\ref{W0}]
of the usual Laplacian (coupled to the bundle $E$).
Using (\ref{fun}) again, we get
	\be\label{w0}
\hf\Del_h=\hf\Del+\hf|\dr h|^2-\Lam_3(\form).
	\ee

When the bundle $E$ is flat, the only difference of $\Boxb_v$, $\Boxb_{\ii v}$
and $\hf\Del_h$ are in the terms $\Lam_a(\form)$ ($a=1,2,3$).
In this case, deformations break the $SU(2)$ symmetry of $\Boxb$ to $U(1)$,
while the $\frac{\pi}{2}$ rotations in $SU(2)$ interchanges the three
operators $\Boxb_v$, $\Boxb_{\ii v}$ and $\hf\Del_h$.

Finally we come to the relation of $\Boxb_h$ and $\Boxb_v$.

\begin{prop}\label{conj}
	\be
\con^{-1}\Boxb_h\con=\Boxb_v-\Lam_3(\ii\Om)-\Lam_1(\ii\Om)-\ii r_V+\ii\Lhat.
	\ee
\end{prop}

\proof{Using (\ref{mag}), (\ref{conbox}) and (\ref{conlam}), we get
	\bea
S_2(\al)^{-1}\Boxb_hS_2(\al)\eq
\Boxb-(1-\cos\al)\Lam_3(\ii\Om)+\sin\al\Lam_1(\ii\Om)+\hf|\dr h|^2	\nno
\sep{-}\cos\al\Lam_3(\form)-\sin\al\Lam_1(\form)-\ii r_V+\ii\Lhat.
	\eea
Set $\al=-\frac{\pi}{2}$.}

	\sect{Localization to the fixed-point set}

\begin{defn}
For $u>0$, $T\ge 0$, let $P\uT(x,x')$ ($x,x'\in M$) be the smooth kernel
associated to the operator $\exp(-u^2\Boxb_{Th/u}+\ii uT\Lhat)$
calculated with respect to the Riemannian volume element $\dr v_M$ of $M$.
\end{defn}

So for $x\in M$, $P\uT(x,x)\in{\rm End}(\Om^{0,*}(M,E))|_x$.
Moreover, $\e{\th\Lhat}P\uT(\e{-\ii\th}x,x)\in{\rm End}(\Om^{0,*}(M,E))|_x$.

\begin{prop}\label{Loc}
Take $\al>0$. There exist $c>0$, $C>0$ such that for all $x\in M$
with $d(x,F)\ge\al$, $\e{\ii\th}\in S^1$, and all $u\in(0,1]$, we have
	\be\label{loc}
|P\uTu(\e{-\ii\th}x,x)|\le c\,\e{-C/u^2}.
	\ee
\end{prop}

\proof{We use the techniques (and the notations) of [\ref{BL}].
Consider $i\colon F\to M$ as an embedding of compact complex manifolds.
Let $\eta=i^*E$ and $\xi_k=\Lam^kT^{*(1,0)}M\otimes E$ ($k=0,\cdots,n$).
Then
	\be\label{res}
(\xi,i_v)\colon 0\to\xi_n\to\xi_{n-1}\to\cdots\to\xi_0
	\ee
is a holomorphic chain complex of vector bundles on $M$.
Since $F$ is discrete, (\ref{res}), together with the restriction map
$\xi_0|_F\to\eta$, is a resolution of the sheaf $i_*\OO_F(\eta)$.
The elliptic operator considered in [\ref{BL}] is
	\be\label{bl}
u^2\Boxb_{Tv/u}=(uD^M+T\Vhat)^2=u^2(D^M)^2+uT\{D^M,\Vhat\}+T^2\Vhat^2
	\ee
acting on $\Om^{*,0}(M)\hat{\otimes}\Om^{0,*}(M,E)=\Om^{*,*}(M,E)$,
where $D^M=\pdrb_v+\pdrb^*_v$, $\Vhat=i_v+i^*_v$.
Particularly important is that the operator $\{D^M,\Vhat\}$ is of
order zero, hence $uT\{D^M,\Vhat\}$ is uniformly bounded for $u\in(0,1]$,
$T\in[0,1/u]$.
We now extend the domain of our operator $\Boxb_h$ from (the 
$L^2$-completion of) $\Om^{0,*}(M,E)$ to (that of) $\Om^{*,*}(M,E)$.
Since the operator preserves the bi-grading of $\Om^{*,*}(M,E)$, it
suffices to prove (\ref{loc}) for the heat kernel with the extended domain.
Using Proposition \ref{conj}, we have
	\be\label{blr}
\con^{-1}(u^2\Boxb_{Th/u}-\ii uT\Lhat)\con=u^2\Boxb_{Tv/u}-r\uT.
	\ee
Here $r\uT=u\Lam_3(\ii\Om)+u\Lam_1(\ii\Om)+uT\ii r_V$ is also uniformly
bounded for $u\in(0,1]$, $T\in[0,1/u]$.
The operator on the right hand side of (\ref{blr}) has the same heat kernel
$P\uT$ up to a conjugation by $\con$.
Therefore the proof of [\ref{BL}, Proposition 11.10] implies that there exist
a sufficiently small $b>0$ (determined by the injectivity radius of $M$),
and $c_1>0$, $C_1>0$ such that for all $x_0\in M$, $u\in(0,1]$, $T\in[0,1/u]$,
$x\in B(x_0,b/2)$, we have
	\be
|(P\uT-P\uT^{x_0})(x,x)|\le c_1\e{-C_1/u^2}.
	\ee
Here $P\uT^{x_0}$ is the smooth heat kernel of the same operator
with Dirichlet conditions on $\pdr B^M(x_0,b)$.
Hence
	\be\label{loc1}
|(P\uTu-P\uTu^x)(x,x)|\le c_1\e{-C_1/u^2}
	\ee
for all $x\in M$ and $T\ge 0$.
(The condition $T\le 1$ can be lifted by a scaling argument.)
Since $\Vhat$ is invertible on $M\m F$, by the proof of
[\ref{BL}, Proposition 12.1], for any $\al>0$ there exist $c_2,C_2,C'_2>0$
such that
	\be\label{loc2}
|P\uTu^x(x,x)|\le\frac{c_2}{u^{2n}}\e{-C_2T^2/u^2+C'_2T}
	\ee
for any $x\in M$ with $d(x,F)\ge\al$.
(\ref{loc1}) and (\ref{loc2}) imply that for some $c,C>0$,
	\be\label{loc3}
|P\uTu(x,x)|\le c\,\e{-C/u^2}.
	\ee
Formula (\ref{loc}) follows from [\ref{Bjdg}, equation (12.7)]:
	\be\label{cauchy}
|P\uTu(\e{-\ii\th}x,x)|\le|P\uTu(\e{-\ii\th}x,\e{-\ii\th}x)|^{\inv{2}}
|P\uTu(x,x)|^{\inv{2}}
	\ee
and from the $S^1$-invariance of $|P\uTu(x,x)|$.}

Clearly, Proposition~\ref{Loc} is valid without the assumption that the
fixed-point set $F$ is discrete; in general $F$ is a symplectic,
hence K\"ahler submanifold of $M$.
The result can also be proved using the method of [\ref{R}].
The proof here is similar to that of [\ref{BZ}, Theorem 3.11] except
that, without the $\ZZ_2$ symmetry there, we do not get a vanishing result
in Proposition \ref{Lim} below.

\begin{defn}
Let $R_p(\th)$ be the isotropy representation of $\e{\ii\th}\in S^1$
on $T_pM$ and let $Z=(z^1,\cdots,z^n)$ be the linear complex coordinates
on $T_pM$ such that the action of $R_p(\th)$ is
	\be
R_p(\th)(z^1,\cdots,z^n)=(\e{\ii\lam^p_1\th}z^1,\cdots,\e{\ii\lam^p_n\th}z^n).
	\ee
For $T\ge0$, set
	\be
\BTp=\inv{2}\Del^p+\inv{2}T^2\sum^n_{k=1}|\lpk|^2|z^k|^2
	+T\sum^n_{k=1}\ii\lpk(e^ke^\kb-i_ki_\kb)
	\ee
and
	\be
\CTp=\inv{2}\Del^p+\inv{2}T^2\sum^n_{k=1}|\lpk|^2|z^k|^2
	-\inv{2}T\sum^n_{k=1}\lpk([e^k,i_k]+[e^\kb,i_\kb]),
	\ee
where $\Del^p$ is the (positive) flat Laplacian on $T_pM$.
\end{defn}

It is easy to see that the $SU(2)$ group elements $S_a(\al)$ ($a=1,2,3$)
act on $\Om^{*,*}(T_pM)$ and that
	\be
\con^{-1}\CTp\con=\BTp.
	\ee
This can be used to recover [\ref{Bams}, Theorem 1.6] from [\ref{W0}].
Moreover, if $\QTp$ is the heat kernel associated to $\exp(-\CTp)$,
then $\con^{-1}\QTp\con$ is that of $\exp(-\BTp)$.

\begin{prop}\label{Lim}
For $T>0$, $\th\in\RE$,
	\bea
\vc\lim_{u\to0}\Tr_{\Omk(M,E)}\exp[-u^2\Boxb\Thu+(\thT)\Lhat]	\nno
\eq\sump E_p(\thT)\Tr_{\Omk(T_pM)}[R_p(\th)\exp(-\CTp)].	\label{lim}
	\eea
Moreover, the limit is uniform in $\th$.
\end{prop}

\proof{We recall the notations of [\ref{BL}, \S 11-12].
Fix a small $\eps>0$.
For $p\in F$, the ball $B^{T_pM}(0,\eps)\subset T_pM$ is identified with
the ball $B^M(0,\eps)\subset M$ by the exponential map.
Let $k'(Z)=\det(\dr_Z\exp)$, $Z\in T_pM$, be the Jacobian.
Then $\dr v_{T_pM}(Z)=k'(Z)\dr v_M(Z)$ and $k(0)=1$.
We also identify $T_ZM$, $E_Z$ with $T_pM$, $E_p$, respectively,
by the parallel transports along the geodesic connecting $p$ and $Z$.
The operators $D^M$ and $\Vhat$ now act on smooth sections of
$\Lam^*(T_pM)\otimes E_p$ over $B^{T_pM}(0,\eps)$.
The setup here is simpler than that in [\ref{BL}, \ref{Bjdg}] because $F$ is
discrete and because of the resolution (\ref{res}) we choose.
(Using the notations in [\ref{BL}, \S 8.f], here $\xi^+=0$ and $\xi^-=\xi$.)
Following [\ref{BL}, \S 11.h-i and \S 12.d-e], we define
	\be\label{L1}
L^{1,p}\uT=u^2(1-\rho^2(Z))\frac{\Del^p}{2}+\rho^2(Z)(u^2\Boxb_{Tv/u}-r\uT),
	\ee
where $\rho(Z)=\rho(|Z|)$ is a smooth function such that $\rho(Z)=1$ if
$|Z|\le\frac{\eps}{4}$ and $\rho(Z)=0$ if $|Z|\ge\frac{\eps}{2}$, and
	\be\label{L3}
L^{3,p}\uT=F^{-1}_uL^{1,p}\uT F_u,
	\ee
where $F_u$ is a rescaling: $F_uh(Z)=h(Z/u)$.
Let $P\uT^{1,p}(Z,Z')$, $P\uT^{3,p}(Z,Z')$ be the smooth heat kernel
associated to the operators $\exp(-L^{1,p}\uT)$, $\exp(-L^{3,p}\uT)$,
respectively, calculated in the volume element $\dr v_{T_pM}$.
Clearly
	\be\label{scale}
u^{2n}P\uT^{1,p}(uZ,uZ')=P\uT^{3,p}(Z,Z').
	\ee
The only term in $L^{3,p}\uT$ that did not appear in
[\ref{BL}, equation (11.60)] is $-\rho^2(uZ)r\uT(uZ)$.
It is easy to see that for $u\in(0,1]$, $T\in[1,1/u]$, the operator
$1_{u|Z|\le\eps/2}r\uT(uZ)$ is uniformly bounded with respect to the
norm $|\cdot|_{u,T,0,0}$ in [\ref{BL}, Definition 11.23].
This, together with [\ref{BL}, Proposition 11.24], is enough to establish
the results in [\ref{BL}, Theorem 11.26] (in the special case of $Z_0=0$) for
$L\uT^{3,p}$.
We can then proceed as the proof of [\ref{BL}, Theorem 11.31] and obtain the
analog of [\ref{BL}, Theorem 12.14] on the uniform estimates of $P\uTu^{3,p}$.
In particular, for any $m\in\NN$, there exists $c>0$ such that if $u\in(0,1]$,
then
	\be
|P\uTu^{3,p}(Z,Z)|\le\frac{c}{(1+|Z|)^m}
	\ee
for $|Z|\le\frac{\eps}{8u}$.
Using (\ref{scale}) and the analog of (\ref{cauchy}), we get
	\be\label{bound}
u^{2n}|P\uTu^{1,p}(u\Rinv Z,uZ)|\le\frac{c}{(1+|Z|)^m}.
	\ee
Next, from (\ref{L1}) and (\ref{L3}), we get
	\be\label{L3u}
L\uTu^{3,p}=\hf u^2(1-\rho^2(uZ))\Del^p+\rho^2(uZ)(D^M)^2
+\rho^2(uZ)(T\{D^M,\Vhat\}+u^{-2}T^2\Vhat^2(uZ)-r\uTu(uZ)).
	\ee
It is easy to see that as $u\to0$, $r\uTu(uZ)\to\ii Tr_V(p)$; the rest of the
terms in (\ref{L3u}) tends to $\BTp$ by [\ref{BL}, Propositions 12.10, 12.12].
Hence
	\be
L\uTu^{3,p}\to\BTp-\ii Tr_V(p),\quad\mbox{as}\,\,u\to0.
	\ee
Proceed as the proof of [\ref{BL}, Theorem 12.16] (with the simplification
$L_{u,1}=L\uTu^{3,p}$ and $L_{u,2}, L_{u,3}, L_{u,4}=0$) and
as [\ref{BL}, \S 12.i], we conclude that
	\be
\con^{-1}P\uTu^{3,p}\con\to\con^{-1}\QTp\con\otimes\e{\ii Tr_V(p)},
\quad\mbox{as}\,\,u\to0
	\ee
in the sense of distributions on $T_pM\times T_pM$.
By the uniform estimates on $P\uTu^{3,p}$,
	\be\label{limP3}
\e{\th\Lhat}P\uTu^{3,p}(R^{-1}(\th)Z,Z)\to 
R_p(\th)\QTp(\Rinv Z,Z)\otimes\e{(\th+\ii T)r_V(p)},\quad\mbox{as}\,\,u\to0
	\ee
uniformly in $\th$ and in $Z$ belonging to any compact set in $T_pM$.
Using (\ref{scale}) and taking the (local) trace over anti-holomorphic
forms only, we get
	\bea
\vc\lim_{u\to0}u^{2n}
\tr_{\Om^{0,k}_p\otimes E_p}[\e{\th\Lhat}P\uTu^{1,p}(uR^{-1}(\th)Z,uZ)]  \nno
\eq E_p(\thT)\tr_{\Om^{0,k}_p}[R_p(\th)\QTp(\Rinv Z,Z)].	\label{limP1}
	\eea
The arguments leading to (\ref{loc1}) imply (see [\ref{BL}, \S 12.d] and
[\ref{Bjdg}, \S 12.d]) that there are $c_0,C_0>0$ such that for all
$u\in(0,1]$ and $Z\in T_pM$ with $|Z|\le\frac{\eps}{8}$, we have
	\be
|P\uTu((p,\Rinv Z),(p,Z))k'(Z)-P\uTu^{1,p}(\Rinv Z,Z)|\le c_0\e{-C_0/u^2}.
	\ee
Therefore in (\ref{limP1}) and (\ref{bound}), $P\uTu^{1,p}(u\Rinv Z,uZ)$
can be replaced by $P\uTu((p,u\Rinv Z),(p,uZ))k'(uZ)$
for $|Z|\le\frac{\eps}{8u}$.
By the dominated convergence theorem (as in [\ref{BL}, Remark 12.5], but
adapted to take into account uniform convergence), we get
	\bea
\vc\lim_{u\to0}\int_{B^M(F,\eps/8)}\tr_{\Om^{0,k}_x\otimes E_x}
[\e{\th\Lhat}P\uTu(\e{-\ii\th}x,x)]\dr v_M(x)			\nno
\eq\sump E_p(\thT)\int_{T_pM}
\tr_{\Om^{0,k}_p}[R_p(\th)\QTp(\Rinv Z,Z)]\dr v_{T_pM}(Z)
	\eea
uniformly in $\th$.
By Proposition \ref{Loc}, we can replace the domain of the integration
on the left hand side by $M$ and thus the proposition follows.}

\vspace{-2ex}

\begin{defn}
Let $q(u,\th)=\sum_{m\in\zz}q_m(u)\e{\ii m\th}$ be a family of formal
characters of $S^1$ parameterized by $u\in\RE$ and
let $q(\th)=\sum_{m\in\zz}q_m\e{\ii m\th}\in\RE((\e{\ii\th}))$.
We say that $\lim_{u\to u_0}q(u,\th)=q(\th)$ in $\RE((\e{\ii\th}))$
if for all $m\in\ZZ$, $\lim_{u\to u_0}q_m(u)=q_m$.
\end{defn}

\begin{cor}\label{equiv}
For $T>0$, the limit (\ref{lim}) holds in $\RE((\e{\ii\th}))$.
\end{cor}

\proof{From Proposition \ref{Lim}, we know that as $u\to0$,
	\be
\Tr_{\Omk(M,E)}\exp[-u^2\Boxb\Thu+(\thT)\Lhat]-
\sump E_p(\thT)\Tr_{\Omk(T_pM)}[R_p(\th)\exp(-\CTp)]\to0
	\ee
uniformly in $\th$, and hence in $L^2(S^1)$ as well.
This implies that all the Fourier coefficients of the left hand side
tend to $0$.
The result follows.}

	\newpage
	\sect{Proof of the theorem}

As explained in the introduction, the heat kernel proof of equivariant
Morse-type inequalities is based on the following

\begin{lemma}\label{MS}
For $u>0$, $T>0$, we have
	\be\label{ms}
\sumk t^k\Tr_{\Omk(M,E)}\exp(-u^2\Boxb\Thu+\th\Lhat)=
\sumk t^k H^k(\th)+(1+t)Q_{u,T}(\th,t)
	\ee
in $\RE((\e{\ii\th}))[t]$ for some $Q_{u,T}(\th,t)\ge0$.
\end{lemma}

\proof{Recall that $\Boxb\Thu=\{\pdrb\Thu,\pdrb^*\Thu\}$.
Since $\pdrb\Thu$ and $\pdrb$ differ by an $S^1$-invariant conjugation,
their cohomologies are isomorphic as representations of $S^1$.
Using the ($S^1$-equivariant) Hodge decomposition, we get
	\bea\label{hodge}
\Tr_{\Omk(M,E)}\exp(-u^2\Boxb\Thu+\th\Lhat)\eq H^k(\th)
+\Tr_{\pdrb^*\Thu\Om^{0,k+1}(M,E)}\exp(-u^2\pdrb^*\Thu\pdrb\Thu+\th\Lhat)  \nno
\sep{+}\Tr_{\pdrb\Thu\Om^{0,k-1}(M,E)}\exp(-u^2\pdrb\Thu\pdrb^*\Thu+\th\Lhat)
	\eea
as formal characters of $S^1$.
Notice that the spectrum of the operator $\pdrb^*\Thu\pdrb\Thu$ on
(the closure of) $\pdrb^*\Thu\Om^{0,k+1}(M,E)$ is identical to that of
$\pdrb\Thu\pdrb^*\Thu$ on (the closure of) $\pdrb\Thu\Om^{0,k-1}(M,E)$.
Since the $S^1$-action commutes with all the operators, we obtain
	\bea
\vc\Tr_{\pdrb^*\Thu\Om^{0,k+1}(M,E)}\exp(-u^2\pdrb^*\Thu\pdrb\Thu+\th\Lhat)\nno
\eq\Tr_{\pdrb\Thu\Om^{0,k}(M,E)}\exp(-u^2\pdrb\Thu\pdrb^*\Thu+\th\Lhat)\ge0
	\label{pos}
	\eea
in $\RE((\e{\ii\th}))$.
We denote either of the expressions in (\ref{pos}) by $Q^k_{u,T}(\th)$.
Summing over $k=0,\cdots,n$ in (\ref{hodge}), we obtain (\ref{ms}) with
$Q_{u,T}(\th,t)=\sumk Q^k_{u,T}(\th)t^k\ge0$.}

We now take the limit $u\to0$. To use Proposition \ref{Lim} or Corollary
\ref{equiv}, we need the following result on the equivariant heat kernel
of the anti-holomorphic sector of the (supersymmetric) harmonic oscillator.

\begin{lemma}\label{cal}
For $T>0$,
	\be\label{2nd}
\Tr_{\Omk(T_pM)}[R_p(\th)\exp(-\CTp)]=\sum_{I\subset\{1,\cdots,n\},|I|=k}
\frac{\e{-T\sumk|\lpk|-T\sum_{k\not\in I}\lpk+\ii\th\sum_{k\in I}\lpk}}
{\prodk[(1-\e{-(T-\ii\th)|\lpk|})(1-\e{-(T+\ii\th)|\lpk|})]}.
	\ee
\end{lemma}

\proof{The operator $\CTp$ acting on $\Om^{0,*}(T_pM)$ splits
$S^1$-equivariantly to $n$ copies of
	\be\label{2d}
\CT=\hf\Del+\hf T^2|\lam|^2|z|^2
-\hf T\lam(-1+[\dr\zb\wedge,i_{\pdr/\pdr\zb}]),
\quad \lam\in\ZZ\m\{0\}
	\ee
acting on $\Om^{0,*}(\CO)$.
Here $S^1$ act on $\CO$ by $R(\th)=\e{\ii\lam\th}$ and hence on
$\Om^{0,*}(\CO)$ as well.
(\ref{2d}) is the sum of the Hamiltonian for the two dimensional harmonic
oscillator
	\be
\HT=\hf\Del+\hf T^2|\lam|^2|z|^2
	\ee
and a bounded operator of order zero.
The smooth heat kernel associated to the operator $\exp(-\HT)$ acting on
$\Om^0(\CO)$ is given by Mehler's formula
	\be
K_{T^2}(z,z')=\frac{T|\lam|}{2\pi\sinh T|\lam|}
\exp\left[-T|\lam|\left(\frac{|z|^2+|z'|^2}{2\tanh T|\lam|}
-\frac{{\rm Re}(\zb z')}{\sinh T|\lam|}\right)\right].
	\ee
Therefore
	\bea
\Tr_{\Om^0(\co)}[R(\th)\exp(-\HT)]
\eq\int_\co\dr^2z\,K_{T^2}(\e{-\ii\lam\th}z,z)			\nno
\eq\frac{\e{-T|\lam|}}{(1-\e{-(T-\ii\th)|\lam|})(1-\e{-(T+\ii\th)|\lam|})}.
	\eea
The bounded 0-th order operator in (\ref{2d}) takes values $T\lam$ and $0$,
respectively, on 0- and 1-forms.
Furthermore, the $S^1$-action $R(\th)$ picks up a phase
$\e{\ii\lam\th}$ on $\dr\zb$.
Therefore
	\be
\Tr_{\Omk(\co)}[R(\th)\exp(-\CT)]=\two{\frac{\e{-T|\lam|-T\lam}}
{(1-\e{-(T-\ii\th)|\lam|})(1-\e{-(T+\ii\th)|\lam|})}}{k=0}
{\frac{\e{-T|\lam|+\ii\lam\th}}
{(1-\e{-(T-\ii\th)|\lam|})(1-\e{-(T+\ii\th)|\lam|})}}{k=1.}
	\ee
Returning to the problem on $T_pM$, 
for $I=\{i_1,\cdots,i_k\}\subset\{1,\cdots,n\}$,
set $\dr\zb^I=\dr\zb^{i_1}\wedge\cdots\wedge\dr\zb^{i_k}$.
Then we have
	\be\label{I}
\Tr_{\Om^0(T_pM)\dr\zb^I}[R_p(\th)\exp(-\CTp)]=
\frac{\e{-T\sumk|\lpk|-T\sum_{k\not\in I}\lpk+\ii\th\sum_{k\in I}\lpk}}
{\prodk[(1-\e{-(T-\ii\th)|\lpk|})(1-\e{-(T+\ii\th)|\lpk|})]}.
	\ee
The trace on $\Omk(T_pM)$ is the sum of (\ref{I}) over $I$ with $|I|=k$.}

(\ref{2nd}) should be interpreted, after a Taylor expansion on the right
hand side, as an equality of formal characters of $S^1$.
\vspace{2ex}

\noindent {\bf Proof of formula (\ref{strong-}):} $\quad$
In (\ref{ms}) we replace $\th$ formally by $\thT$ and still regard it
as an equality of formal series in $\e{\ii\th}$.
Since as $u\to0$ the limit of the left hand side exists in $\RE((\e{\ii\th}))$
(Corollary \ref{equiv}) and since $H^k(\th)$ is independent of $u$,
we conclude that $\lim_{u\to0}Q_{u,T}(\thT,t)=Q_T(\thT,t)$ also exists
and that $Q_T(\th,t)\ge0$.
Therefore
	\be
\sump E_p(\thT)\Tr_{\Omk(T_pM)}[R_p(\th)\exp(-\CT)]
=\sumk t^kH^k(\thT)+(1+t)Q_T(\thT,t).
	\ee
Using Lemma \ref{cal} and changing $\thT$ back to $\th$, we get
	\be
\sum_{p\in F,I\subset\{1,\cdots,n\}}t^{|I|}E_p(\th)
\frac{\e{-T\sumk|\lpk|-T\sum_{k\not\in I}\lpk+(T+\ii\th)\sum_{k\in I}\lpk}}
{\prodk[(1-\e{\ii|\lpk|\th})(1-\e{-|\lpk|(2T+\ii\th)})]}	\nno
=\sumk t^kH^k(\th)+(1+t)Q_T(\th,t).
	\ee
Finally, as $T\to+\infty$, the limit of each summand on the left hand side is
$0$ except when the pair $(p,I)$ satisfies $I=\set{k}{\lpk>0}$, in which case
the limit is 
$t^{n-n_p}E_p(\th)\e{\ii\sum_{\lpk>0}|\lpk|\th}/\prodk(1-\e{\ii|\lpk|\th})$.
Consequently, $\lim_{T\to+\infty}Q_T(\th,t)=Q^-(\th,t)$ exists as well,
and $Q^-(\th,t)\ge0$.
Formula (\ref{strong-}) is proved. $\hfill\Box$

\newpage	

\noindent {\bf Acknowledgement.} 
We are grateful to Weiping Zhang for bringing to our attention and kindly
explaining the works in [2-6], and for critical comments on the manuscript.
We also thank Sixia Yu for discussions.
S.\ W.\ thanks the hospitality of the University of Adelaide,
where this work started.
The work of S.\ W.\ is partially supported by NSF grant DMS-93-05578
at Columbia University and by ICTP.

\bigskip

        \newcommand{\athr}[2]{{#1}.\ {#2}}
        \newcommand{\au}[2]{\athr{{#1}}{{#2}},}
        \newcommand{\an}[2]{\athr{{#1}}{{#2}} and}
        \newcommand{\jr}[6]{{#1}, {\it {#2}} {#3}\ ({#4}) {#5}-{#6}}
        \newcommand{\pr}[3]{{#1}, {#2} ({#3})}
        \newcommand{\bk}[4]{{\it {#1}}, {#2}, ({#3}, {#4})}
        \newcommand{\cf}[8]{{\it {#1}}, {#2}, {#5},
                 {#6}, ({#7}, {#8}), pp.\ {#3}-{#4}}
        \vspace{5ex}	
        \begin{flushleft}
{\bf References}
        \end{flushleft}
{\small
        \begin{enumerate}
        
        \item\label{AB}
        \an{M.\ F}{Atiyah} \au{R}{Bott}
        \jr{A Lefschetz fixed point formula for elliptic complexes, Part I}
        {Ann. Math.}{86}{1967}{374}{407};
        \jr{Part II}{Ann. Math.}{87}{1968}{451}{491}

	\item\label{Bi} 
	\au{J.-M}{Bismut} 
	\jr{The Witten complex and the degenerate Morse inequalities} 
	{J. Diff. Geom.}{23}{1986}{207}{240}

	\item\label{Bams}
	\au{J.-M}{Bismut}
	\jr{Koszul complexes, harmonic oscillators, and the Todd class}
	{J. Amer. Math. Soc.}{3}{1990}{159}{256}

%
	\item\label{Bjdg} 
	\au{J.-M}{Bismut}
	\jr{Equivariant immersions and Quillen metrics}
	{J. Diff. Geom.}{41}{1995}{53}{157}

	\item\label{BL} 
	\an{J.-M}{Bismut} \au{G}{Lebeau}
	\jr{Complex immersions and Quillen metrics}
	{Inst. Hautes Etudes Sci. Publ. Math.}{74}{1991}{1}{297}

	\item\label{BZ}
	\an{J.-M}{Bismut} \au{W}{Zhang}
	\jr{Real embeddings and eta invariants}
	{Math. Ann.}{295}{1993}{661}{684}

	\item\label{D0}
	\au{J.\,P}{Demailly}
	\cf{Sur l'identit\'e de Bochner-Kodaira-Nakano en g\'eom\'etrie
	hermitienne}{S\'eminaire d'analyse, Lecture Notes in Math., vol.\ 1198}
	{88}{97}{ed.\ \athr{P}{Lelong et.\ al.}}{Springer-Verlag}
	{Berlin and New York}{1986}

	\item\label{D}
	\au{J.\,P}{Demailly}
	\jr{Champs magn\'etiques et in\'egalit\'es de Morse pour la 
	$d''$-cohomologie}{Ann. Inst. Fourier (Grenoble)}{35}{1985}{185}{229};
	\cf{Holomorphic Morse inequalities}
	{Several complex variables and complex geometry, Part 2,
	Proc. Symposia Pure Math., vol. 52}{93}{114}
	{eds.\ \athr{E}{Bedford} et.\ al.}
	{Amer. Math. Soc.}{Providence, RI}{1991}

	\item\label{Fr}
	\au{T}{Frankel}
	\jr{Fixed points and torsion on \ka manifold}
	{Ann. Math.}{70}{1959}{1}{8}

	\item\label{GH}
	\an{P}{Griffiths} \au{J}{Harris}
	\bk{Principles of Algebraic Geometry}
	{John Wiley \& Sons, Inc.}{New York}{1978}, \S 0.7

	\item\label{L}
	\au{K}{Liu}
	\jr{Holomorphic equivariant cohomology}
	{Math. Ann.}{303}{1995}{125}{148}

	\item\label{R} 
	\au{J}{Roe}  
	\bk{Elliptic operators, topology and asymptotic methods,
	{\rm Pitman Res.\ Notes in Math.\ 179}}
	{Longman Scientific \& Technical}{Essex, England}{1988}, \S 12

	\item\label{W0}
	\au{E}{Witten} 
	\jr{Supersymmetry and Morse theory}
	{J. Diff. Geom.}{17}{1982}{661}{692}
	
	\item\label{W} 
	\au{E}{Witten} 
	\cf{Holomorphic Morse inequalities}
	{Algebraic and differential topology, Teubner-Texte Math., 70}
	{318}{333}{ed.\ G.\ Rassias}{Teubner}{Leipzig}{1984}

	\item\label{Wu}
	\au{S}{Wu}
	\pr{Equivariant holomorphic Morse inequalities II: torus and
	non-Abelian group actions}
	{MSRI preprint No.\ 1996-013, {\tt dg-ga/9602008}}{1996}
	
        \end{enumerate}
	\end{document}